\documentclass[12pt]{article}

\usepackage{graphicx}

\usepackage{amssymb}
\usepackage{amsmath}

\begin{document}

\newtheorem{theorem}{Theorem}[section]
\newtheorem{lemma}{Lemma}[section]
\newtheorem{corollary}{Corollary}[section]
\newtheorem{claim}{Claim}[section]
\newtheorem{proposition}{Proposition}[section]
\newtheorem{definition}{Definition}[section]
\newtheorem{fact}{Fact}[section]
\newtheorem{example}{Example}[section]

\newcommand{\quod}{\hfill $\blacksquare$ \bigbreak}
\newcommand{\reals}{I\!\!R}
\newcommand{\property}{I\!\!P}
\newcommand{\np}{\mbox{{\sc NP}}}
\newcommand{\sing}{\mbox{{\sc Sing}}}
\newcommand{\con}{\mbox{{\sc Con}}}
\newcommand{\prob}{\mbox{Prob}}
\newcommand{\atm}{\mbox{{\sc ATM}}}
\newcommand{\hopn}{\hop_{\cN}}
\newcommand{\atmn}{\atm_{\cN}}
\newcommand{\cA}{{\cal A}}
\newcommand{\cO}{{\cal O}}
\newcommand{\cP}{{\cal P}}
\newcommand{\cC}{{\cal C}}
\newcommand{\C}{{\cal C}}
\newcommand{\cR}{{\cal R}}
\newcommand{\cB}{{\cal B}}
\newcommand{\cG}{{\cal G}}
\newcommand{\cN}{{\cal N}}
\newcommand{\cU}{{\cal U}}
\newcommand{\cF}{{\cal F}}
\newcommand{\cT}{{\cal T}}
\newcommand{\hx}{\hat{x}}
\newcommand{\cS}{{\cal S}}

\newcommand{\MS}{\mathcal{S}}
\newcommand{\CR}[1]{\mbox{cr}\left( #1, n \right)}

\newcommand{\eps}{{\epsilon}}
\newcommand{\la}{{\lambda}}
\newcommand{\al}{{\alpha}}
\newcommand{\qed}{\hfill $\square$ \smallbreak}

\newenvironment{proof}{\noindent{\bf Proof:}}{\qed}

\def\lalto{\left \lceil}
\def\ralto{\right \rceil}
\def\lbasso{\left \lfloor}
\def\rbasso{\right \rfloor}
\def\D{{\Delta}}
\def\qed{\hfill$\Box$}

\baselineskip    0.28in
\parskip         0.1in
\parindent       0.0in

\bibliographystyle{plain}
\title{{\bf Why Do We Believe Theorems? }}
\author{
Andrzej Pelc \footnotemark[1] }
\date{ }
\maketitle
\def\thefootnote{\fnsymbol{footnote}}

\footnotetext[1]{ \noindent
D\'epartement d'informatique, Universit\'e du Qu\'ebec en Outaouais, Gatineau,
Qu\'ebec J8X 3X7, Canada.
E-mail: {\tt pelc@uqo.ca}  
}

\begin{abstract}
We investigate the reasons of having confidence in mathematical theorems.
The formalist point of view maintains that formal derivations underlying 
proofs, although usually not carried out in practice, contribute to this confidence.
Opposing this opinion, the main claim of the present paper is 
that such a gain of confidence obtained from any link
between proofs and formal derivations is, even in principle, impossible in the present state of knowledge.
Our argument is based on considerations concerning length of formal derivations.
We also discuss psychological and social factors that contribute to building belief in theorems.

This is a significantly extended version of a paper with the same title published in:\\
{\em Philosophia Mathematica} III 17 (2009), 84-94.\\
and reprinted in:\\
The Best Writing on Mathematics 2010,\\ Mircea Pitici, Ed., Princeton University Press, 2010.

\vspace*{7cm}

\end{abstract}

\pagebreak

\section{Introduction}

Most of the common mathematical practice deals with proofs of theorems. As 
authors, mathematicians\footnotemark[1]
\footnotetext[1]{The term ``mathematician'' is used here in a large sense and 
encompasses all scientists
adopting the methodology of proofs taken from ``mainstream'' mathematics. 
Hence, apart from, say,
algebraists or topologists, it includes, e.g., mathematical logicians and 
theoretical computer scientists but excludes users of mathematics, such as 
engineers or physicists, interested in mathematical results 
but not in the way they are obtained. To justify this classification through a 
personal example, the present author,  formerly a set theorist and currently 
working in theoretical computer science, has not felt  significant 
methodological change when changing his subject of interest: the (high-level) 
methodology remained the same
-- proving theorems then, proving theorems now, only subjects of these 
theorems changed, although quite dramatically, from large cardinals to 
distributed algorithms.}
invent proofs and try to write them down rigorously; as readers, they try to 
verify and understand proofs of other mathematicians; as referees and 
journal editors, they assess the interest and value of proofs, and as 
teachers, they explain proofs to novices.

Why are mathematicians so concerned with proofs? Or, reformulated in a more direct way:
Why do we prove theorems? This question serves as the title of the paper by Rav [1999].
At first, one may be tempted to give a deceptively simple answer: we prove theorems to 
convince ourselves and others that they are true. While this is often, indeed, the direct reason
for proving a theorem (very often a conjecture is first stated by someone and then the same person
or some other mathematician proves it to give it a status of a theorem, or disproves it, to give a status
of a theorem to its negation), the convincing power is far from being the unique role of proofs
in the mathematical practice. Often the new ideas and techniques conveyed by a proof are much more important
than the theorem for which the proof was originally invented. 
Rav [1999] (page 20) formulates it succinctly:

\begin{quotation}
Proofs are for the mathematician what experimental procedures are for the experimental scientist: 
in studying them one learns of new ideas, new concepts, new strategies -- devices which can 
be assimilated for one's own research and be further developed.
\end{quotation}

and illustrates this epistemic function of proofs by well chosen examples.

Nevertheless, the role of proofs as means of convincing the mathematical community about 
the validity of theorems is very important. While proofs can also serve other purposes,
only proofs can directly serve this purpose. Since the subject of the present study is establishing
reasons for believing theorems, rather than discussing epistemic value of proofs,
our perspective on proofs is much more restrictive than that of Rav [1999].
He dealt with the question ``Why do we prove theorems?'' and so was concerned with all aspects of proofs,
both as ways to convince that the proved theorem is correct, and as repositories of mathematical knowledge.
Our question (the title of this paper is deliberately modeled on the title of Rav [1999]) is
``Why do we believe theorems?''. Thus we are only interested in the ``convincing'' role of proofs.
However, since in many
cases a particular mathematician does not have access to a proof of a particular result,
and still has to form a belief about it, we will also consider other factors that contribute to such beliefs.
As mentioned above, only proofs can {\em directly} serve the purpose of convincing mathematicians about
the validity of theorems but, as discussed in Section 5, there are also indirect ways to form such convictions, 
on which mathematicians often have to rely.
 
Before announcing our position concerning reasons of believing theorems, we need to explain two crucial terms
that will be used throughout the paper: {\em proof} and {\em derivation}. In choosing these particular
terms, we follow Rav [1999]. (The expressions {\em informal proof} instead of {\em proof} and
{\em formal proof} instead of {\em derivation}  could be also used, as they convey the sense we want to 
give to those terms, but we find Rav's terminology more convenient.) By {\em proofs} we will mean the arguments 
used in mathematical practice in order to justify correctness of theorems. Since this
simply reports the meaning adopted by the mathematical community, there is no formal definition of this term.
On the other hand, the term {\em derivation} is used in its formal, logical sense. 
Recall that a formalized language is
described first. It has an appropriate alphabet (containing, among others, logical symbols). 
Terms and formulae of this language are defined 
as specific strings of symbols of this alphabet. Then a formalized theory $\cal T$ is defined in this language
by specifying a -- possibly infinite -- recursive set\footnotemark[1]
\footnotetext[1]{This is a technical term whose informal meaning is that there exists a mechanical method
to verify if a given sentence is a member of the set or not.} 
of formulae called {\em axioms} that form the basis of the theory. A finite set of transformations called 
{\em inference rules} is specified. These transformations permit to obtain a new formula from a finite set 
of formulae. One example of an inference rule is {\em modus ponens}. Finally, a (linear) {\em derivation} 
is a finite sequence of formulae such that every term of this sequence is either an axiom or is obtained from
a set of earlier terms of this sequence by applying one of the inference rules. A formula of the above
formalized language is called a {\em theorem} of the theory $\cal T$, if it is the last formula in some derivation.
The crucial characteristic of this formal definition is that, due to the fact that the set of axioms is recursive
(although possibly infinite), there is a mechanical way to verify if a given sequence of formulae of the 
formalized language is a derivation or not. (This should not be confused with the possibility of a mechanical
method to verify if a given formula is a theorem of $\cal T$ or not. The well-known undecidability theorem states
that the latter cannot be verified by a mechanical method for sufficiently strong theories.) This feature
of derivations is the reason of their important role in the construction of formalized theories. Once a 
(alleged) derivation is presented, there is a mechanical way to verify its correctness, and in the 
case of a positive verdict, to assert that the last formula of the derivation is indeed a theorem of the theory $\cal T$.

\section{Statement of position}

Since the formulation of Hilbert's program, most of the discussions concerning reasons for confidence
in mathematical theorems are conducted along the axis between formalists and their opponents. In a nutshell, 
the first group asserts that {\em derivations} underlying proofs actually carried out 
in mathematical practice are the reason to believe those theorems, while the second group dismisses
or significantly weakens the role of those derivations in building confidence in theorems, and points out
to other features of proofs, such as the analysis of the meaning of mathematical notions, as well as to 
social factors, as responsible for the beliefs in mathematical statements. Hence, while the radical interpretation
of Hilbert's program, as expressed in the famous Problem Number 2\footnotemark[1] 
\footnotetext{Hilbert [1900]:
I wish to designate the following as the most important among the numerous questions which 
can be asked with respect to the axioms: 
{\em To prove that they are not contradictory, that is, that a finite number of logical steps based upon 
them can never lead to a contradictory result.} (Italics in original)}
stated in Hilbert's address in Paris in 1900, has been shown unfeasible by the later negative results of G\"odel,
a more moderate interpretation, trying to harvest confidence gains from some link between proofs and derivations, 
is still present and lively in the discussions of the topic.

This controversy between formalists and their opponents is well exemplified by a discussion between 
Jody Azzouni and Yehuda Rav in a sequence of papers Rav [1999], Azzouni [2004], and Rav [2007]. In this
discussion Azzouni presents a moderate formalist point of view which he calls the
{\em derivation-indicator} view, while Rav opposes his arguments. Azzouni's view is well summarized
in the abstract of Azzouni [2004], where he says:
\begin{quotation}
A version of Formalism is vindicated: Ordinary mathematical proofs indicate 
(one or another) mechanically checkable derivation of theorems from the assumptions 
those ordinary mathematical proofs presuppose. The indicator view explains why 
mathematicians agree so readily on results established by proofs in ordinary 
language that are (palpably) not mechanically checkable.
\end{quotation}
The present paper is not an exception from the rule that reasons of mathematical beliefs are discussed
along the controversy between formalists and non-formalists. We choose Azzouni [2004] as a protagonist
of the formalist point of view because his position is moderate\footnotemark[2]: 
\footnotetext[2]{For example, it is further said in the abstract of Azzouni [2004]:
``Mechanically checkable derivations in this way structure ordinary mathematical practice 
without its being the case that ordinary mathematical proofs can be 'reduced to' such derivations.''}
we want to show that even such a moderate
position is impossible to maintain. However, our opposition to the views expressed by Azzouni [2004]
is presented from a standpoint much different from that of Rav [2007]. In fact, we will argue that the critique of 
Azzouni as proposed by Rav [2007] is too harsh in some aspects and too mild in others. To summarize our position,
we will argue that: 
\begin{enumerate}
\item
If reasonably short derivations could be carried out for all mathematical theorems then 
significant gains of confidence in these theorems could be obtained.
\item
The premise of the above statement is not known to be true in the present state of knowledge.
\item 
No link between derivations and proofs can contribute to increasing 
confidence in theorems in the present state of knowledge.
\item
Confidence in theorems may be a dynamic attitude and is obtained through a social process.
\end{enumerate}
Hence we think that Rav's critique of the formalist position is too strong when he says (Rav [2007], page 306,
italics in the original):
\begin{quotation}
Let F(*) denote the resulting formalization of the proof of Theorem (*) as a derivation, 
obtainable in the way as has just been indicated calling upon, as must have been noted, 
a logician's know-how. One can write now a computer program for checking mechanically 
whether the derivation is valid, or, alternatively, code the proof in the AUTOMATH 
language of [de Bruijn, 1970] and use its computer implementation for the mechanical checking. 
Suppose now that the  verdict of the computerized proof-checker is `{\em not valid}'. 
Do we conclude that Theorem (*) is false, or do we conclude that somewhere along the way 
in the formalization or computer program there is a `bug' ? (Rhetorical question.) 
Next suppose that the computer verdict says `{\em valid}'. Do we now have more confidence in the `truth'  
(or validity) of Theorem (*)? Are we now reassured that the (informal) proof with which we have 
started is correct, or have we thereby just been reassured that the whole exercise of 
formalization and mechanical checking have been carried out correctly? (Once more a rhetorical question.)
\end{quotation}
and continues in a footnote:
\begin{quotation}
As a matter of fact, when a human agent has transcribed an informal proof as a derivation 
in a suitably chosen formal system, there is a conceivable loss in the reliability in 
the formalized derivation by comparison to the original informal proof, since the correct 
reading of the informal proof by the human agent and its resulting transcription cannot be 
mechanically checked (at the pain of an infinite regress!).
\end{quotation}

Indeed, while we agree that a {\em hypothetical} negative computer verdict would be inconclusive
(it would amount to the statement: ``you missed this time, try again''), 
we will argue that a {\em hypothetical} positive computer verdict would increase our confidence,
not in the informal proof that has been successfully formalized as a derivation, but in the theorem itself.
However, the key point in our argument, that positions the present author firmly in the antiformalist camp,
is that such a computer verdict is impossible. (Note the italics in the word ``hypothetical'', used above twice.) 

On the other hand, when Rav says (Rav [2007], page 294, italics in the original):
\begin{quotation}
Consequently, when it comes to the nature of the {\em logical} justification of mathematical 
arguments in proofs, with Azzouni putting his faith in {\em formal derivations}, 
even if just indicated, and further maintaining that on this basis proofs are recognized 
by mathematicians as valid, as they see them hence being {\em mechanically checkable}, 
against such formalist-mechanist claims I have objections to voice, both on technical 
and on historical grounds. As opposed to various formalist views, I hold that mathematical 
proofs are {\em cemented} via arguments based on the {\em meaning} of the mathematical terms that 
occur in them, which by their very conceptual nature cannot be captured by formal calculi.      
\end{quotation}
his critique is too mild. Indeed, we will argue that not only do mathematical proofs contain ingredients
that cannot be captured by formal calculi but, in fact, no link between those proofs and hypothetical
derivations underlying them can be established, in case of many theorems. 
Thus, while Rav argues that there is added epistemic value
in proofs with respect to derivations, we make a much stronger claim: in the case of many theorems, 
derivations cannot contribute any epistemic value at all.

\section{The convincing power of short derivations}

In this section we justify our first claim: if short derivations were possible to carry out for 
all mathematical theorems, then 
significant gains of confidence in these theorems could be obtained. Later on we will argue that
the premise of this claim is not true, so, from the point of view of formal logic, 
the above claim seems to have no interest,
as every sentence is implied by a false statement. However, this claim is quite important to explain our position
regarding the convincing power of derivations. It shows in particular that the possibly prohibitive length
of derivations is the {\em only} reason why confidence in theorems cannot be gained from a (hypothetical)
link between proofs and derivations. This is in sharp contrast with the above quoted opinion of Rav [2007]
that even if such derivations could be carried out and (using his terminology) a positive ``computer verdict'' 
concerning the correctness of the derivation 
could be given, this would not increase our confidence, and that such 
confidence would even be decreased in the course of formalization, as compared to the confidence drawn from the
original (informal) proof.

In order to justify our claim, we propose the following thought experiment. 
For the purpose of the experiment we fix an axiomatized theory, for example
Zermelo-Fraenkel set theory with the Axiom of Choice (ZFC), in which derivations will
be carried out. Moreover we assume (this is for the purpose of the experiment only, we would certainly
not maintain such an absurd position) that
all theorems published in mathematical papers have short derivations in ZFC, 
of 1000 steps at most.
Of course, any previously published result can be used in such derivations and counts as one step.
Hence in order to obtain a ``pure'' derivation of a theorem added late in the process,
i.e., a derivation that would not contain
references to previously published results, the 1000-step derivation of it would have to be ``unfolded''
by inserting derivations of all previous theorems that it uses, these derivations should be again ``unfolded'',
and so on. The final ``pure'' derivation (not using any previous theorems) could of course be quite long, much longer
than 1000 steps, but not arbitrarily long: any step of it had to be recorded some time earlier in the process,
by some previous author.
This crucial fact will be used later in our negative argument.

Coming back to our thought experiment, 
we further assume that such derivations can be easily, although tediously, obtained given a standard
well-written (informal) proof. Hence we suppose a situation similar to that in computer science, where
an average programmer can write a program in a computer language, given a (thoroughly written) informal
algorithm.

Suppose also that given this hypothetical situation, all mathematical journals introduced a new requirement
at the time of submitting a paper for publication. As an appendix to the paper itself, the author is asked 
to attach derivations (in ZFC) of all results of the paper. Before the paper is sent to the referees,
a technical assistant to the editor feeds all appended derivations to a computer program that checks if 
a given derivation is correct. If all derivations pass the test, the paper goes through the normal refereeing
process. If not, the paper is sent back to the author with the indication which steps of which derivations
failed the test and with a suggestion of correcting errors and resubmitting the paper. The checking program 
itself is in public domain, is the same for all journals, and is itself thoroughly checked by the mathematical
and computer science community. (Such a program {\em verifying} derivations would be itself quite short
and easy to check for correctness, in the informal sense. Writing such a program is an easy programming
exercise).

To continue with our experiment, if the paper is finally accepted (according to the usual refereeing
criteria, no changes here), at the time of its publishing, every theorem of the paper receives a unique
Correctness Identification Number (CIN) published together with the theorem. When another (or the same) 
author wants to use this result in a subsequent paper, the CIN of this theorem has to be quoted in the 
attached derivation. All derivations are deposited in a central Archive of Derivations 
(possibly maintained by Google?). When the computer program checking derivations encounters 
a step quoting an existing result, it verifies its CIN in the central archive and checks if the quoted formulation
in the currently verified step matches exactly the formulation of the theorem with this CIN, 
kept in the archive. In this way the new derivation is incorporated in the universal bank of derivations,
and can be added to the archive when the new paper is published.

Our thought experiment is completed, it is time to discuss its results. Before arguing in what way such a hypothetical
system would increase our confidence in published mathematical theorems, we should point out what it would
{\em not} accomplish. (We are far from suggesting that even such a hypothetical scenario would be a miraculous
remedy for all doubts concerning mathematical truth.) First it should be clear that no epistemic gains
in the sense of increased mathematical knowledge could be obtained from such a system: we would not learn new 
ideas or techniques, we would not gain any insight in the importance of the results of a given paper.
Further, the system would not even necessarily increase the confidence in the correctness of the published {\em proof}
(only, as we will argue, in the proved {\em result}). Indeed, since the relation between the (informal) proof
and the accompanying derivation would not be externally controlled, it is theoretically possible (although highly unlikely)
that a flawed proof would somehow lead to a correct derivation. (Such a situation, although possible, is unlikely,
because the flaw in the proof would be probably discovered and corrected during 
the production of the correct derivation.) 

It should be also clear that a rejection of the derivation by the checking program does not mean that the result
is false, or even that the informal proof on which the derivation was based is incorrect. It only points to errors
in the derivation
that could be hopefully corrected by the author before resubmission. This is similar to the usual debugging
process in the production of computer programs. Finally, such a system would certainly not provide any
deciding procedure for either the truth or provability of mathematical statements. 
(It is well known that such a procedure cannot exist.) The modest goal of the system would just be to verify
if a derivation of a result (provided by the author) is correct and to maintain the (gigantic) graph of such 
derivations with links pointing from a given derivation to previously accepted results used in it.
  
The only potential gain from the above described hypothetical verification system could be the gain 
of confidence in published theorems. Let us analyze it more closely. What kind of confidence
should we have in a result, seeing that it got published in the hypothetical system and thus got 
a CIN certificate? Looking at the above described process, the reader will hopefully agree that:

- if we believe the axioms of ZFC and the axioms of logic,\\
- if we consider that the inference rules used in derivations, such as modus ponens, do not decrease
our belief in the infered statement as compared to the premises,\\
- if we believe that the checking program is correct,\\ 
- if we believe that the hardware of the computer on which the checking program runs is correct,\\
- if we believe that the technical assistant is honest (does not, for example, type a false correctness
certificate),

then we should believe the theorem whose derivation obtained a positive verdict. This is not to say that
we should believe the correctness of the proof on which the derivation is based: as justly pointed
out by Rav [2007] (see the quotation above):
``the correct reading of the informal proof by the human agent and its resulting transcription cannot be 
mechanically checked''. All that we claim is that, under the above five assumptions, we should believe
{\em the theorem itself}. This is a modest claim: let us examine in what way the belief would actually be
increased with respect to the current (real) situation.

There are five conditions on our newly acquired (hypothetical) belief: the first two concern philosophical
belief in axioms and logical rules, the third concerns the belief in the correctness of a very simple
computer program, the fourth concerns the belief in the soundness of technology, and 
the fifth concerns the belief in
the honesty of a particular human. Hence our belief is clearly not absolute. In order to analyze how strong 
it should be, we need to examine all the conditions, as the chain is only as strong as its weakest link.

Concerning the first two conditions it is enough to say that we cannot expect the belief in a 
(possibly complicated) theorem to be stronger than the belief in axioms of ZFC and the axioms of logic, that are simple, 
intuitively obvious statements, and stronger than the belief in the soundness of rules like modus ponens.
It seems that denying the belief in the above precludes the belief in any mathematical statements
and thus renders our discussion uninteresting. (Here we put aside the issue of doubt in some particular, more
problematic axioms of ZFC, like the Axiom of Choice, not accepted by all mathematicians in its strongest form.
For most results in mathematics, a quite weak form of the Axiom of Choice is sufficient,
a form accepted by most mathematicians, and we might modify 
the ZFC system in our experiment to such a weaker system.)

The third condition may be considered a weak point in our network of beliefs: the checking program may be 
obviously incorrect. However, as mentioned before, the checking program would be very simple, and -- being
in public domain -- could be independently checked by thousands of people. Thus, although {\em theoretically}
it would be prone to errors, as any result of human activity, the {\em practical} confidence in it (after years of testing
and using) would be very high.

As for the fourth condition, it is enough to say that hardware faults in computers are extremely rarely
the cause of erroneous results; human factors are much more often to blame in this respect. Another argument is that,
if we believe (as we probably all do) in fault-free functioning of computer hardware in monitoring 
life-saving medical equipment and flight control systems, we should not be more suspicious in the case
of computer derivation checking.

Finally, the technical assistant honesty issue: a delicate point, since, as we know, there is no lack
of dishonest individuals in the world. Here the argument is similar as above: we trust the honesty 
of pharmacists to give us the correct medicine and the honesty of aircraft pilots when they claim 
that they are not under the
influence of drugs or alcohol when flying the plane. 
Their dishonesty could cost lives. Should we be less trusting of the honesty of
editor assistants?

We have hopefully established that the hypothetical system verifying derivations would {\em practically} 
confer quite a high degree of confidence in the theorems in mathematical literature. Would this confidence
be higher than in the actual current situation? In other words, could we then talk about a {\em gain}
in confidence? This author feels that the answer to this question is ``yes''. Let us consider the following
situation in which we believe that our confidence in a theorem would indeed increase, if we could 
use the hypothetical system. 
Suppose that we want to use a little-known result published in an obscure journal. The proof is long,
difficult to understand, and far from our domain of expertise. In fact, we are not interested in the proof at all,
we need the result itself, as a ``black box'' in our argument. We do not have the time, desire, 
or ability to go through the
proof in order to check that it is correct. 
In this case, a CIN of the obscure theorem (obtained via the above described process) would definitely 
increase our confidence by substantially lowering the risk that the result we are about to use is wrong.

Unfortunately, as already stressed several times above, the derivation verifying system is purely hypothetical.
We will argue in the next section that such a system could not be implemented. Indeed, the crucial assumption
of our thought experiment, that short derivations are possible to carry out for 
all (published) mathematical theorems, is blatantly false.

\section{Impossibility of gaining confidence from derivations}

In this section we present the main negative claim of this paper. We will argue that,
in the present state of knowledge, it is
impossible to gain any confidence in most mathematical theorems from hypothetical derivations
underlying their proofs. Hence we will argue against the statement of Azzouni [2004] who says
(page 83):

\begin{quotation}
it's derivations, derivations in one or another algorithmic system, 
which underlie what's characteristic of mathematical practice: 
in particular, the social conformity of mathematicians with respect to 
whether one or another proof is or isn't (should be, or shouldn't be) convincing.
\end{quotation}

Recall that we concluded from our thought experiment in Section 3 that if  
reasonably short derivations of mathematical results could be explicitly written, then
they would contribute to a significant gain of confidence in the theorem under scrutiny.
Hence it is not surprising that our impossibility argument concerning such a gain in reality
is based on length considerations concerning such hypothetical derivations.

It should be stressed that our impossibility claim is quite strong. We want to argue that
not only mathematicians do not in practice use derivations to get or increase confidence
in their results, but that in the present state of knowledge 
it is {\em theoretically impossible} to achieve such a gain of confidence,
in the case of most interesting mathematical theorems.

Since we want to argue against the above quoted statement of Azzouni, we should try to understand
what he means by saying that derivations ``underlie what's characteristic of mathematical practice''.
First we should ask: what exactly is an algorithmic system and derivations 
in which algorithmic system? Azzouni [2004] provides a kind of answer to the first question,
saying in the footnote to the above quote:
``An algorithmic system is one where the recognition procedure for proofs 
is mechanically implementable. By no means are algorithmic systems restricted 
to language-based axiom systems.'' He is more vague as far as the choice of a  
specific system is concerned (page 93):
\begin{quotation}
...the picture doesn't require mathematicians, in any case, 
when studying a subject-matter, to remain within the confines of a single 
(algorithmic) system indeed, if anything, mathematicians are required to 
transcend such systems by embedding them in larger ones. 
The derivation indicated (by the application of new tools to a given subject matter) 
can be a derivation of the weaker system the mathematician started with, 
or it can be a derivation of a stronger system (some of) the terms of which 
are taken to pick out the same items supposedly referred to in the weaker system.
\end{quotation}
Nevertheless, for a given theorem in mind, we have to focus on a particular 
``algorithmic system'' (while agreeing that this system could be modified several 
times in the course of developing the proof of the theorem). 
Derivations in an axiomatized theory, such as
Zermelo-Fraenkel set theory with the Axiom of Choice (ZFC) are in principle mechanically checkable
(when written correctly), hence we will use the system ZFC as an example of an algorithmic system
in which derivations are ``indicated'' by usual proofs, to use Azzouni's terminology. ZFC is a good
example because many mathematicians believe that most of the mathematical lore could be 
{\em in principle} formalized in this system. Actually, defining other mathematical notions 
(such as functions, relations, algebraic systems) in terms of sets is a standard mathematical
practice which can be considered as a step in the direction of such a formalization. (Further steps are
usually not made). 

Next we should see what Azzouni means by saying that proofs ``indicate'' derivations. This point is 
somewhat obscure, as a precise definition of this relation between proofs and derivations is never given
in his paper. To be sure, Azzouni [2004] is well aware of the fact that ``mechanically recognizable derivations ...
(generally) aren't ones (that can be) exhibited {\em in practice}'' (page 95, italics in the original).
However, the very fact of his stressing that derivations are not exhibited {\em in practice} points to 
the important distinction between practical considerations and a {\em theoretical} accessibility of
such derivations. This is also implied by his statement ``it's {\em derivation} which provides the skeleton
for (the flesh of) {\em proof}'' (page 95, italics in the original). In order to provide
such a skeleton, a derivation should be {\em at least theoretically} accessible for scrutiny,
and hence it should be {\em at least theoretically} possible to {\em record and verify} it. 
Since Azzouni stresses mechanical checkability as an important feature of derivations, 
it is thus fair to assume that,
whatever could be the precise meaning of the phrase ``proofs indicate derivations'', 
it should be {\em theoretically possible} to mechanically check such indicated derivations.

Now comes a crucial remark that forms the basis of our impossibility argument. Since Azzouni (and other formalists)
claim that confidence concerning theorems is somehow linked to derivations underlying proofs,
the burden of justifying that such (indicated) derivations are (at least theoretically) possible to be recorded and 
mechanically checked is on {\em their} side, in every case when such a gain of confidence is claimed.
A caveat is in order here. Formalists do not need to actually {\em exhibit} an indicated derivation,
but they should justify that such a derivation is theoretically possible to be recorded and 
mechanically checked.  
This theoretical possibility is a {\em necessary} condition for such a gain of confidence and hence should be justified
by anybody claiming the gain. There is an important asymmetry with respect to the ``burden of justification''.
In order to counter this formalist claim for a given theorem, there is no need to show that such a
mechanically checkable derivation is 
theoretically impossible; it is enough to show that the formalist side has not justified the theoretical
possibility of mechanical checking of such a derivation. We will argue that this is the case, and hence that 
the claimed confidence gain cannot be harvested.   
 
Before proceeding with our impossibility argument, a disclaimer is in order. Our argument refers
to ``complex'', ``deep'' theorems, as opposed to simple corollaries from definitions, such as
the unicity of the neutral element in groups, or simple geometric theorems, e.g.,  
that the sum of angles
of a triangle totals $\pi$. For the latter two examples it is not hard to imagine a derivation
in ZFC; it is probably even possible (although very likely a terribly boring task) 
to write it explicitly and thus, according to our thought experiment from Section 3, which 
should be possible to actually carry out in such simple cases, to gain confidence in those results.
(As a counter-argument we might notice that confidence in these simple results is extremely high
anyway, so the additional potential gain would be small, if not non-existent.) Of course,
the above notions of ``complex'', ``deep'' theorems are very fuzzy, 
which is the reason to put them in quotation marks. 
However, hardly any mathematician will disagree that the Fermat Last Theorem proved by Wiles [1995]
is both complex and deep (no need of quotation marks in this case). Hence we will use this theorem,
refered to as $FLT$, as an example in our impossibility argument.

We should also stress that the choice of ZFC as the underlying axiomatic system 
in which derivations are indicated by usual proofs is only given as an example. It could be
replaced by any normally used system such as Peano Arithmetic, Kelley-Morse theory of classes, etc.,
whose axioms and inference rules people find intuitively obvious. 
Of course one might create an artificial system, e.g.,  containing FLT as an axiom
(in which the derivation of FLT would have length one),
but such manipulations cannot possibly increase confidence in this theorem.

The starting point of our argument is the following controversy between Rav [1999] and Azzouni [2004].
Rav [1999] (page 14) writes:

\begin{quotation}
In reading a paper or monograph it often happens as everyone knows too well that one arrives 
at an impasse, not seeing why a certain claim $B$ is to follow from claim $A$, as its author affirms. 
Let us symbolise the author's claim by `$A \rightarrow B$' . (The arrow functions iconically: 
there is an informal logical path from $A$ to $B$. It does not denote formal implication.) 
Thus, in trying to understand the author's claim, one picks up paper and pencil and 
tries to fill in the gaps. After some reflection on the background theory, 
the meaning of the terms and using one's general knowledge of the topic, 
including eventually some symbol manipulation, one sees a path from $A$ to 
$A_1$, from $A_1$ to $A_2$, $\dots$, and finally from $A_n$ to $B$. This analysis can be 
written schematically as follows: 
$$A \rightarrow A_1,A_1 \rightarrow A_2,\dots,A_n \rightarrow B.$$ 
Explaining the structure of the argument to a student or non-specialist, 
the other may still fail to see why, for instance, $A_1$ ought to follow from $A$. 
So again we interpolate $A \rightarrow A', A' \rightarrow A_1$ . But the process of interpolations for 
a given claim has no theoretical upper bound. In other words, how far has one 
to analyse a claim of the form `from property $A$, $B$ follows'
before assenting to it depends on the agent. There is no theoretical reason to 
warrant the belief that one ought to arrive at an atomic claim $C \rightarrow D$ 
which does not allow or necessitate any further justifying steps between $C$ and $D$. 
This is one of the reasons for considering proofs as infinitary objects. 
Both Brouwer and Zermelo, each for different reasons, stressed the infinitary character of proofs.
\end{quotation}

and Azzouni [2004] (page 97) responds:

\begin{quotation}
What's going on? Why should anyone think that a finitary piece of mathematical reasoning, 
a step in a proof, say, corresponds to something infinitary (if, that is, 
we attempt to translate it into a derivation)?
\end{quotation}

While we agree with the main thought of Rav [1999], it seems that the 
use of the word ``infinitary'' may have been exaggerated. It is this
word and not the argument itself that seems to have provoked the controversy. 
What we consider as the crux of Rav's reasoning is that there is no 
clearly defined {\em upper bound} on the number of justifying steps. If the
argument is presented in this way, there is no need to invoke any infinitary 
character of the proof, that raised Azzouni's objections. Our impossibility argument
will be based on this important distinction. (In the sequel we never make any claims 
about an infinitary character of proofs or derivations. Similarly as Azzouni, we believe that
these are finite objects.) 

For any theorem $T$ whose derivation (or many derivations) in ZFC are indicated 
(according to Azzouni's terminology), by its ``ordinary'' published proof, let ${\cal L} (T)$ denote the
length of the shortest possible derivation (not only among those derivations indicated by the 
particular considered proof, but among all possible derivations of $T$ in ZFC). 
Such derivations exist since, by assumption, they are indicated by the ordinary proof. Hence
the integer ${\cal L} (T)$ is well defined. We conjecture that in the case
of most ``complex'', ``deep'' theorems $T$, and in particular in the case of $FLT$, it is impossible to 
provide (and justify) {\em any} upper bound on ${\cal L} (T)$. We agree that this is a bold statement but we would
like to challenge a skeptical reader to provide (and justify!) any upper bound on ${\cal L} (FLT)$. 
Should it be a million? a quadrillion? $10^{10^{10}}$? What would be the justification of the choice of any such number?
Actually, of any number at all? To be sure, we have not {\em proved} that it is impossible to give such an estimate,
we only conjectured it. Nevertheless, and this time it is not a conjecture but a statement of a fact: 
nobody (including formalists claiming confidence gains from derivations underlying proofs) 
has ever given any such estimate. This simple observation is the first crucial claim in our argument.  
Our conjecture is stronger and says that such an estimate is 
impossible to establish, but it will not be used in the argument: due to the asymmetry with respect 
to the ``burden of proof'', the fact that such an estimate has not been given will be enough for our reasoning. 

The next step of our argument requires the definition of a very large integer number $M$, 
an {\em extremely} large number, indeed. We are not concerned with determining the size of this number.
The only thing that matters for our purposes is that no derivation of length larger than $M$
could ever be actually written down or verified mechanically, even in theory. 
We will proceed with a series of estimates leading to the
definition of $M$. First consider the period of time known as the {\em Planck time} (cf. Halliday et al. [1996]).
This interval of time, call it $t_P$, of length of order $10^{-43}$ seconds, is a lower bound on the duration 
of any observable
physical event, and hence it is also a lower bound on the time needed for any conceivable processing device
to record or verify one step of a derivation. Next, switching to the macroscale, consider a future state of
the universe precluding any information processing activity. A good candidate for such a state would be
the heat death of the universe, i.e., when the universe has reached maximum entropy. Let $\tau$ denote an upper bound
(in seconds) on the time till this state of the universe is reached. 
Let $Q=\tau /t_P$. Hence $Q$ is an upper bound on the number of 
steps of a derivation that can ever be recorded or verified by a single processing device. It is of course  
possible that many processing devices cooperate in verifying some derivation in parallel, each processing unit working on
a different part of the derivation. Hence we need a third estimate, an upper bound on the number of 
elementary particles in the observable universe. Call this number $P$. Clearly, any processing device must be composed of
at least one particle, hence the number of such conceivable devices is at most $P$. Finally let $M= P \cdot Q$.
It follows that the maximum number of derivation steps that could ever be recorded or verified is bounded by $M$.
In other words, a derivation of length larger than $M$ could never be verified, even if the entire observable
universe were converted into a machine totally devoted to the verification of this single derivation.

It should be noted that all our estimates leading to the definition of $M$ are grossly exaggerated:
the time of recording or verifying one step of a derivation is much larger than $t_P$, any computing activity in
the universe would stop sooner than after $\tau$ seconds, and the number of processing devices in a hypothetical verifying
machine is much smaller than $P$. As a consequence, our claim concerning the number $M$ could be even made about a much smaller
number. However, we chose to use the above exaggerated estimates because, for the purpose of our argument, we only need to
indicate some number bounding the length of a theoretically verifiable derivation, and the number $M$ has the advantage
of being rather simple to define and of not leaving any doubt concerning the validity of our claim about it.       

Consider a derivation of length larger than $M$.
The hypothetical existence of such a derivation of a theorem $T$ could not possibly contribute to our
confidence in $T$ because we could never have any kind of access (even theoretically) to
all terms of such an extremely large sequence of formulae, and hence we could never verify that it is indeed
a correct derivation. Trying to gain
any confidence concerning the validity of theorem $T$ from the existence of such a hypothetical derivation
is (metaphorically speaking) similar to getting insight into a black hole.

Let us now call a theorem $T$ {\em reachable}, if ${\cal L} (T) \leq M$. In other words, a theorem is reachable,
if and only if there exists its derivation in ZFC of length at most $M$. It follows from what was said above that 
reachability of a theorem is a necessary condition for gaining any confidence in it from some hypothetical derivation
of this theorem.   

Now we are ready for the third, final step of our argument. Consider the theorem $FLT$. Since no upper bound
on ${\cal L} (FLT)$ has been justifiably provided, the answer to the question
of whether $FLT$ is reachable is unknown. Given the fact that the number $M$ is so enormous, one would be 
tempted to give the answer ``yes'' to this question, i.e., to establish $M$ as an upper bound on
${\cal L} (FLT)$. However, as observed before, this has never been done 
(and seems impossible to do, although we do not need this stronger conjecture in our argument). 
To avoid misunderstandings: we are not claiming
that $FLT$ is not reachable, we are only arguing that the answer to this question is unknown.

We can now conclude that, since in the case of $FLT$ (and the reasoning would be similar in the case of
many other ``complex'', ``deep'' theorems) we are unable to 
decide -- in the present state of knowledge -- whether $FLT$ is reachable, there
can be no gain of confidence from indicating any derivation of $FLT$. Indeed, as observed at the beginning of our argument,
the side claiming such a gain has the burden of justifying that for some derivation of $FLT$ it is 
{\em theoretically possible} to mechanically check it. This however would imply reachability of $FLT$,
which, as previously argued, is not known to be true. 

We would like to present the following analogy.
A biologist studying animals' behavior or a trainer working with animals considers 
such notions as fear, anger, pain, or sexual attraction,
and tries to explain various actions of animals (attacking, running away, mating, etc.) 
in these terms. 
However, these actions could also be explained at some lower level, the level of biochemical reactions in the
animal's brain, or even at a subatomic level, involving statistical information about motions of particles.
It is clear that these lower levels (especially the subatomic one) are not very useful in explaining and
predicting the animal's behavior: they are too detailed. Likewise, a proof of a theorem should be explained
at a ``macroscopic'' level, involving mathematical objects and principal techniques, rather than at a 
``subatomic'' derivation level. However, we would like to observe that the analogy ends at this point.
Indeed, the subatomic level in the case of living organisms can be (indirectly) observed. 
Hence one may argue that (with suitable hypothetical tools) 
we may be able to gain additional biological knowledge studying such a microscopic level. 
In the case of mathematical theorems, the accessibility of this ``subatomic'' level,
i.e., the level of derivations, is quite different. The core of the problem is their {\em mechanical verifiability}.
As previously argued, in the case of particular theorems, these objects  
may be too large to be written, perceived (in any possible sense), 
let alone to be scrutinized for mechanical verification. So while in biology it is unlikely but theoretically possible
that the subatomic level can provide additional insight into animals' behavior,
such a possibility (even theoretical) may not exist in mathematics in the case of many theorems, for
a simple but crucial reason: the size of objects at the ``subatomic '' level 
(i.e., at the level of derivations) may be too large to be examined with the aim of verification 
by a human mind, or by any technological device.

Finally, it should be stressed that, as previously mentioned, 
our impossibility argument concerns the present state of knowledge.
In other words, we argued that it is {\em now} impossible (even theoretically) 
to gain additional confidence concerning the validity of
theorems such as $FLT$ by ``indicating'' their derivations. Indeed, we based our argument on
the observation that the question of whether
$FLT$ is reachable has not been settled until now. 
This impossibility at present is enough to refute Azzouni's [2004] claim, which
is made about such gains of confidence at present. It is of course hard to predict how the situation could change with 
additional insights in the future. If, for example, short derivations of all published mathematical theorems
were discovered one day (a very unlikely but not completely impossible scenario) then the situation would resemble our 
thought experiment from Section 3 and Azzouni's claim could be substantiated. This, however, does not concern our discussion
whose purpose was refuting the possibility of such confidence gains {\em now}, a claim made by the formalist side.

\section{Psychological and social aspects of believing\\ 
theorems}

It follows from the arguments in Section 4 that derivations cannot increase our confidence
in most important mathematical theorems, in the present state of knowledge. 
Since these formal well-defined objects are of no help,
we have to seek reasons of such belief elsewhere. We will argue that a combination of
psychological and social factors is responsible for our confidence in theorems.
This idea is not new. De Millo, Lipton and Perlis [1979] (page 272) say:
`` ... insofar as it is successful, mathematics is a social, informal, intuitive, organic, 
human process, a community project''. 
Thomas [1990] (page 80) is even more explicit:
\begin{quotation}
A public [mathematical] statement is regarded as true when the mathematical 
community can be convinced that it has been properly deduced from what is 
already explicit in the public mathematics.
\end{quotation}

However, we will provide (hopefully) new arguments for the importance of psychological and social factors
in building (or destroying) belief in theorems. Rather than adopting a static approach that
concentrates on the issue of which elements of proofs are likely to convince mathematicians, we will discuss
dynamic aspects of beliefs, and try to answer the question of how and why such beliefs change in time,
both at an individual level and at the level of the mathematical community. We will also argue that confidence 
in a particular theorem (even in a fixed point in time) is not a totally universal attitude:
it may depend on the individual mathematician. 
    
In order to give an example of psychological and social factors involved in 
forming our beliefs
concerning the validity of theorems, we propose our second 
thought experiment.
Consider the following theorems:

$S_{1}$: Your own recent result that you have just written down and submitted 
to a journal
(here we assume that the reader is a mathematician in the large sense, 
mentioned in the introduction).

$S_{2}$: A Ph.D. thesis result by a brilliant young mathematician, published
in the recent issue of a reputable journal.

$S_{3}$:  An insignificant result in graph theory, with a very complicated 
proof, published in the
proceedings of an obscure conference in 1970 \footnotemark[1] 
\footnotetext[1]{For obvious reasons, no example is given, not because of 
lack thereof.} 
(here we assume that the reader is not an expert in graph theory).

$S_{4}$: The Simple Group Classification Theorem (the so called ``Enormous 
Theorem'').

$S_{5}$: The  Pythagorean Theorem.

$S_{6}$: The Riemann Mapping Theorem (here we assume that the reader is not a 
specialist in complex analysis).

Now suppose that somebody proposes you the following bet, concerning each of the 
above statements
separately. If, within a year from the bet, the statement is shown false (in 
a way so convincing that you will agree
with this yourself) then you pay  \$10,000; otherwise you get  \$10,000. On 
which of these statements would you bet (assuming that you do not dismiss the 
very idea of betting for reasons independent of beliefs concerning the 
subject of the bet, e.g., for moral or religious reasons)?
The advantage of setting the experiment in the framework of a bet is
that the will to bet (especially a large sum of money) reflects 
rather well what the parties believe and how strongly, regardless of the 
reason for such beliefs. Moreover, this setting has the advantage of ``discretization''
of the inherently gradual nature of beliefs: while we may ``somewhat'' or ``almost''
believe a theorem, a bet is necessarily a 0-1 decision.  

Now let us perform a reasoning that, we think, reflects to some extent the 
decision process in which an average
mathematician would engage in this situation. At least this is the reasoning 
that this author would make and the 
reader is invited to compare how he or she would react. Let us look at the 
statements one by one.

$S_{1}$ concerns
{\em your own} result and proof. Shouldn't we believe that our own work is 
correct? Of course we should and 
we do, often too strongly...  How often does one get a referee's 
report saying that Theorem 2.3 is not
correct as it is written because some (hopefully insignificant) assumption is 
missing, or the statement is true not
for all topological spaces but for all except the singleton space. Such 
mistakes are easy to correct but, according
to the rules of the bet, would be sufficient to lose it. Unfortunately, 
referees sometimes also spot more serious mistakes in places completely 
unsuspected by the inventor of the proof. Hence, this author would, 
reluctantly, refuse the bet in the case of  $S_{1}$.
 
Should we bet on $S_{2}$? It looks really good: at least two independent 
experts (we are talking about a reputable journal) agreed that the proof is correct.
However, with some hesitation, this author would likely again
refrain from the bet. After all,  \$10,000  is a lot of money for a university professor, 
and this author has seen a few false results even
in reputable journals. What is lacking in the case of  $S_{2}$ is the 
opportunity for more members of the mathematical community to look at the 
proof, possibly use its ideas and confront it with background knowledge in 
the domain: this is the way errors are often found. Such an opportunity 
requires some time after publication of the proof.  

In the case of $S_{3}$, the main risk lies in the possibility that 
nobody except the author has ever carefully
looked at the proof. Although almost forty years elapsed since its 
publication, the insignificance of the result combined with the length and 
complications of the proof and the fact that the paper is ``buried'' in an 
obscure place,  make this quite plausible. Moreover, conference proceedings 
of this type are sometimes not refereed or refereed
only superficially, so there is a good chance that nobody except the author 
has seen the paper also prior to publication. Hence the danger here is even 
larger than in the case of $S_{2}$, in spite of the age of the result. 
Too much of a risk, at least for this author. 

In the case of $S_{4}$ the 
decision is really difficult. This is a very well-known and thoroughly 
discussed result and (parts of various versions of) the proof have been seen 
by many mathematicians. The proof of the
Simple Group Classification Theorem  is spread over tens of thousands of 
pages in about 500 papers, some of them unpublished. Many mathematicians 
expressed serious doubts about the correctness of the proof at various stages 
of its development, serious gaps were found and subsequently corrected (cf. 
Aschbacher [2004]). Some 
experts still doubt that the current version of the proof is correct. 
Needless to say, personal verification of the entire proof is out of the 
question for almost everyone. Not a betting case for now, according to 
this author.  

Now comes  $S_{5}$: this is a 
moneymaker -- no doubt about jumping on the bet. Why? Is it {\em absolutely} 
impossible for the Pythagorean Theorem to be wrong? Well, not absolutely...
It is {\em theoretically} possible that all people writing, reading and 
teaching various (simple) proofs of this theorem over the ages were 
completely blind to some obvious error in them and that the result itself is 
incorrect because 
of this overlooked error.
There is no theoretical way to exclude such a possibility. Nevertheless, 
everyone would agree that this is so
extremely unlikely that it should not prevent even the 
most cautious person from betting on the theorem. Here we have a
very reassuring combination of evidence: the proof is simple, possible to be 
readily checked by any mathematician (or even non-mathematician), including 
the person about to bet, and historical data show that it was, in fact, 
checked and used by countless people over centuries, without creating a shade of 
doubt.
This is probably the highest degree of confidence (we use this term in an  
informal, everyday sense) that a mathematician can get about the correctness 
of a result.  

The final candidate is $S_{6}$, 
the Riemann Mapping Theorem. This is a very important and well-known result 
in complex analysis, proved by Riemann in 1851 in his Ph.D. thesis. His 
proof, however, contained some hypotheses that were 
subsequently proved to be invalid in some cases. The first correct proof was 
given by Carath\'{e}odory in 1912.
Subsequently  other, simpler proofs were given by other mathematicians and 
the result itself, as well as
techniques of these proofs, were used ever since. Here the betting decision 
is made in a slightly different context
than in the case of 
the Pythagorean Theorem. Most non-specialists in complex analysis have 
probably never 
personally checked the proof, and actually most of them (this author included) 
have never seen it. 
However, the result and its various proofs are very important, old and well 
known, which guarantees
that they have been successfully checked and used by many people. Thus 
(without even spending time
to try to look at any of the proofs by himself) this author would bet on 
$S_{6}$ without hesitation.

Our thought experiment is over: in the case of $S_{5}$ or $S_{6}$ we would 
probably cash a hefty check 
after a year, in other cases we might well regret excessive caution, if the 
young brilliant mathematician
really found a correct proof (or even if the proof was incorrect and the 
result false but no one found any 
mistake in it within a year), if the obscure graph theoretic
result was valid and the proof correct (although long and boring), and if no 
new dramatic negative developments
concerning the Enormous Theorem occurred within a year.

Now let us consider the outcome of the experiment.
Its usefulness lies more in the presented reasoning behind each decision 
(we think that these reasonings reflect, to some extent at least,
the decision making process of other mathematicians in such a hypothetical 
situation), than the decision itself. 
\footnotemark[1]
\footnotetext[1]{Before discussing what follows from our thought experiment 
and to avoid misunderstandings, we should stress what does NOT follow from 
it. Certainly we do not want to imply any kind of formal probability measure 
attached to theorems, in particular connected to some randomized checking of 
proofs, such as mentioned and justly criticized by Rav [2007]. All 
notions of likelihood and risk used here are informal
and purely subjective.}
The above presented arguments indicate that the reasons of a belief in a 
particular mathematical statement are a combination of psychological and 
social factors. We take into account our personal ability to verify the proof 
of a result, the complexity or simplicity of the presented proof which makes 
it more or less likely to hide errors and conceptual mistakes, 
the degree to which the result
is established in the mathematical community (which in turn depends on its 
visibility, importance and the time
since it has been published), the reputation of the venue, which somehow 
indicates the seriousness of the
refereeing process, known controversies surrounding the proof, as in the case 
of the Enormous Theorem,  and many other  factors, not directly spelled out 
above. 

An important consequence of the previously presented 
arguments is that
the belief in a particular result is a dynamic concept: the same result 
becomes more ``worthy of belief'' as time goes by and its proof is reverified 
and used by others, or, conversely, the belief in a result may be shattered,
if a gap in the proof is revealed (some steps in it shown to be doubtful), or 
even worse, if some of the lemmas it uses are shown false. These dynamic 
changes in belief may last for years, as it is the case with the Enormous 
Theorem, depending on the new developments and critiques of the proof. 

It should also be clear that the belief in a theorem depends on the 
qualifications of the person pondering
its validity. Curiously, the strength of this belief is neither an increasing 
nor a decreasing function of these
qualifications. On the one hand, a novice may be blocked  at some step of the 
verification of a proof, 
not seeing why it is really the 
case that ``B clearly follows from A''  as stated in the journal paper, while 
an expert will not have any doubt
about this step because he or she has seen it in other proofs and used it 
many times before. On the other hand,
a novice may ``slide'' over a doubtful argument, not even seeing any dangers, 
where a more experienced reader
immediately ``sniffs'' a real problem.

Hence it is our opinion that the way in which mathematicians gain (or lose) 
belief in particular theorems
does not differ radically from what happens in other sciences. The underlying 
dynamic social process is
similar. When Azzouni [2004] (page 84) says that ``mathematicians are so good at 
agreeing  with one another on whether
some proof convincingly establishes a theorem'', he seems to put 
mathematicians in a completely different
category, than, e.g., specialists in experimental sciences or  engineers.
 We would like to argue that this is not the case: on the one hand 
mathematicians are
not so good at it, on the other hand, others are not so bad. To justify the 
first point let us note three reasons
for disagreement between mathematicians, concerning the validity of a given 
proof: historical reasons, difference of competence, and mistakes in 
arguments. 

The historical aspect has been extensively discussed by Rav [2007]
who quotes, among many other papers, the historical analysis of Grabiner 
[1974].  The latter shows,  that,
using the words of Rav [2007] (page 294) ``rigor, and hence the type of proofs that are 
accepted by the mathematical community, are time-dependent''. 
This time-dependence is not a monotonic function.
On the one hand, some proofs perfectly acceptable in the eighteenth and nineteenth 
centuries, let alone proofs of first results in calculus by Leibniz or 
Newton, would not be accepted in a modern journal because the referee would 
complain that  the author barely conveys some ideas and does not provide any 
``real'' proof of his statements. The rigor requirements in calculus have significantly 
increased since its creation. On the other hand, some type of reasoning rigor 
becomes relaxed with time, as a particular method or technique becomes widely known 
and understood. It is perfectly acceptable to say in a modern set theory paper
``adding $\omega _2$ Cohen reals forces property $P$'', while this would be considered
an obscure statement in the 1960's, right after the creation by P.J. Cohen of the method 
of forcing used to construct models of set theory.  
Hence the agreement as to what is convincing and what is not, does 
not hold across centuries.

The difference of competence and its role in agreeing if some proof is 
convincing,  is well known to all mathematicians. Obviously, because of social 
context, a novice, especially a student, will often not argue
about his or her position very strongly, if there is a controversy about a validity 
of a proof between him or her and some
member of the academic community higher on the hierarchy ladder. The student 
may suppose that the difference
of opinion is due to his or her insufficient knowledge or insight. (This 
timidity may also be a dangerous thing
and should not be encouraged:
often mistakes in proofs are found by young students who try to verify each 
step of the proof because they are 
not sure of anything, and during this slow process find an error, while more 
senior colleagues have just perused the proof and ``convinced themselves'' 
that everything is ``basically correct''.) Nevertheless the
difference of opinions remains a fact, even if it is only temporary. A 
consequence of this is well known 
in the mathematical academic community: a proof presented at an undergraduate 
course must substantially
differ from a proof at a conference. MacKenzie [2001] (page 321) justly writes:

\begin{quotation}
[a] sense of audience is crucial to this form of proving: a sense of what 
listeners or readers will know, and of what
will be enlightening; of what the audience will understand, and what they 
will not; of what needs to be spelt out, and what can be covered by phrases 
such as `it is obvious that...', or `we can show similarly that...', or 
`without loss
of generality, we can assume that...';  and so on. 
\end{quotation} 

A potential lack of this sense of audience may be a reason for disagreement 
over the validity of a result.
In the extreme case when an author of a submitted paper does not succeed in 
writing a proof in a way that convinces the referees, the latter may finally 
reject the paper because they are not sure that the result and 
the proof are correct, and several revisions have not remedied the situation.

The third reason for disagreement are existing mistakes in arguments. The 
example of the Enormous Theorem
shows this issue very explicitly. At what time can we talk about agreement 
that the current proof is correct?
Will a shade of doubt persist in many years to come, at least in the minds of 
some experts?
How can we be sure that the (supposed) belief in correctness will not change 
in the near future? It is tempting to make the ironic statement that at 
any given time the potential agreement can only be about the current version 
of the proof being {\em incorrect} and needing further modifications.

The above arguments show that Azzouni's optimism about the ability of 
mathematicians to reach agreement
which proofs are convincing is somewhat exaggerated. We rather tend to agree with a more
realistic position of Hersh [1979] who says (page 43):
\begin{quotation}
We do not have absolute certainty in mathematics. 
We may have virtual certainty, just as in other areas of life. 
Mathematicians disagree, make mistakes and correct them, 
are uncertain whether a proof is correct or not ...
\end{quotation}  

Now let us briefly 
consider the second
point: are other communities so much worse at reaching similar agreements? 
Let us consider, for example,
the civil engineering community. It does not seem that civil engineers are 
incapable of reaching an agreement
concerning, say, safety standards in bridge construction. If this were the 
case, given the crucial importance
of these standards for human lives, bridges would not be built, engineers 
being eternally occupied 
discussing what standards should be adopted. This does not seem to be the 
case. To give an example in more fundamental sciences: physicists seem to be
in perfect agreement about, say, laws of optics, at least when restricted 
to the everyday environment.   
Examples in other domains
of experimental sciences and technology are not hard to find. 

Hence we would like to argue that the difference between the capability of 
mathematicians to reach agreement
whether a proof convincingly establishes a theorem and a similar capability 
of experts in other domains
of, say, experimental sciences and technology being able to reach agreement 
about the validity of results in these domains does not form a black and 
white picture. When Azzouni [2004] (page 103) says:
\begin{quotation}
It's not that in mathematics we can't get things wrong as traditional philosophers 
have too often suggested; it's that mathematicians have a way of agreeing about proof 
that is virtually unique (agreement of any sort among humans is always a surprise, 
and always has to be explained somehow), and the explanation of this surely is 
that some sort of recognition of mechanical procedures is involved. 
\end{quotation}
we think that Azzouni is wrong both in describing the existing situation and in its explanation. 
There may be a difference of degree in this capability between 
mathematicians and other scientists, but not 
a qualitative difference. In both cases the process is a social one. A given 
community at a given time must be 
convinced, and the way of convincing it is always the same: discussion using 
arguments currently accepted by 
this community. 
\footnotemark[1]
\footnotetext[1]{Cf. Thurston [1994]:
... mathematical truth and reliability come about through the very human process 
of people thinking clearly and sharing ideas, criticizing one another and 
independently checking things out.}

One could argue that pointing to a social process of discussion among experts
as responsible for convincing the mathematical community about the validity of theorems is 
incorrect because belief in theorems and types of acceptable arguments last much longer
than social structures: the Pythagorean Theorem is still believed, although states of 
ancient Greece have fallen long ago. This argument has two flaws. First, the durability of these beliefs
is not always so long. The beliefs in the validity of the Enormous Theorem changed more rapidly than
many modern governments. Also, the social feeling of what is an acceptable argument in a mathematical
journal evolves rather rapidly, as some techniques become widely known and become part of what mathematicians
call {\em folklore}. More importantly, the fact that the belief in the main body of mathematical lore
transcends social {\em structures} does not mean that it is not dependent on a social {\em process}.
This process of scientific discussion, written or oral, is quite similar in Plato's Academy and in a 
modern university department, and its outcomes are, indeed, usually more lasting than the underlying social
structures, e.g., than those institutions themselves or than countries in which they operate.
This does not refute the fact that consensus concerning a result, or lack thereof, are precisely
outcomes of these scientific discussions.   
Of course the types of acceptable arguments vary from domain 
to domain and are time-dependent. The tools used in the discussions are very diversified even inside mathematics 
and include both rigorous parts, such as calculations, 
and more informal ones, such as pictures, diagrams, analogies or even metaphor helping to visualize an argument 
(e.g., a node of a graph is ``killed'', an area of the plane is ``protected'').

The above arguments hopefully convinced the reader that there are no qualitative differences
between the capabilities of mathematicians and of scientists in some other domains to reach an agreement concerning
the validity of an assertion in their respective fields of expertise. However, one other point requires further discussion.
In the case of experimental sciences it may be argued that observation of the external world causes some conformity in
the way specialists in these domains describe and explain phenomena. This raises the question of what is the 
corresponding reason of the (admittedly changing, imperfect and far from absolute) agreement among mathematicians concerning
the validity of a large part of the body of mathematics. (Although we argued that not all mathematicians agree on 
the validity of all theorems at all times, there is no doubt that most mathematicians agree on the validity of
most theorems most of the time.) In the case of mathematics, observation of the external world cannot be given as a reason.

We would like to consider the answer to this important question at two levels: 
first by trying to find the immediate reason for this
general conformity, and then by trying in turn to find the cause of this reason. We think that the immediate reason is this.
At any given point in time the mathematical community accepts
some ``rules of the game'', which are permitted in informal proofs. 
As argued before, these rules change over time: some are relaxed
(when they become folklore), some are made more precise, when the old ones have been shown to contain traps.
However, if we fix a moment in the development of mathematics, these rules are fairly universal 
and respected by most mathematicians,
because most of them are trained on the same proofs,
and learn from the same type of manuals, so they generally accept the same types of arguments.
These rules are not based on derivations, and in fact, were used for centuries 
before the concept of derivation has even been defined,
which happened quite recently (compared to the duration of mathematics).
There was an attempt to formalize and fix
these rules when formal systems and derivations were first defined, but this attempt failed to change
mathematical practice (while being an interesting object of study by itself). 

Now the question arises: why are some rules of reasoning adopted rather than others?
This calls for the second level of our answer. In the chapter ``The Uniqueness of Mathematics as a Social Practice''
of his book Azzouni [2006] writes:
\begin{quotation}
What seems odd about mathematics as a social practice is the presence of substantial conformity on the one hand,
and yet, on the other, the absence of (sometimes brutal) social tools to induce conformity
that routinely appear among us whenever behavior really is socially constrained. Let's call this
``the benign fixation of mathematical practice''.
\end{quotation}
and further gives three reasons for this ``benign fixation'':
\begin{enumerate}
\item
``...a capacity to carry out algorithms, and -- it's important to stress this -- this is a species-wide capacity:
We can carry out algorithms and teach each other to carry (out) specific algorithms {\em in the same way}''
(italics in the original)
\item
``Arithmetic and geometry, in particular, have obviously intended domains of application.
These fixed domains of application to some extent prevent drift in the rules governing terms of mathematics --
in these subjects so applied, anyway -- because successful application makes us loath to change successfully 
applied theorems if that costs us applicability.''
\item
``an at least partially ``hard-wired'' disposition to reason in a particular way'' 
\end{enumerate}

While we agree with Azzouni that the above facts may indeed be true, we do not share his opinion that all of them
are reasons for conformity in mathematics. The capacity to carry out algorithms in the same way does not seem to be an
important reason for this conformity,  due to the fact that only a small part of the mathematical lore has algorithmic
nature. This capacity may well justify the agreement on some arithmetic laws, like commutativity or associativity
of addition, but can be hardly considered a reason for conformity in some abstract parts of mathematics, like
general topology. Likewise, while invoking applications can possibly ``prevent drift'' in rules applied to arithmetic
or geometry, no such applications can be conceived for highly abstract domains like, e.g, the theory of large cardinals,
which are abstract concepts used to differentiate between different sizes of ``very large'' infinite sets.
By contrast, the third reason cited by Azzouni seems more convincing. It is hard not to observe a 
``disposition to reason in a particular way'' widely spread among humans. In other words, some basic rules of
reasoning formalized in logic,\footnotemark[1]
\footnotetext{Here we mean classical logic, since this logic is used in an overwhelming part of mathematics.}
such as modus ponens, or deriving a property of a particular object from the assertion that all objects
have this property, seem to be innate to most human beings. One might even hypothesize that lack of disposition to use
such basic rules during the development of the species would disadvantage an individual in the competition for survival, 
thus causing such dispositions to be ``hard-wired'' in the brain by the evolutionary process. Once a strong majority
adopts such homogeneous rules, these rules are in turn reinforced socially. (It is enough to imagine prospects of
a mathematics student who refuses to use modus ponens or the principle of mathematical induction claiming that these
principles seem unnatural to him or her.) 
If the hypothesis of innate uniform rules of reasoning among humans is true, the homogeneity of these rules 
could play a similar role in preserving (general) conformity in mathematics as observation of the external world 
plays in experimental sciences. 

To summarize, while fundamental reasons for agreement in mathematics and in experimental sciences
may admittedly be different, the social process through which this agreement is reached is similar:
in both cases it is the discussion among experts.     
Both in the case of mathematics and of experimental sciences and 
technological domains, some members
of the community may remain unconvinced, so it is often impossible to reach 
unanimity rather than 
general consensus. Sometimes, virtually all members of the community agree, but 
this is true both in mathematics (in the case of well-known, established
theorems) and in physics for instance (the laws of optics or the non-existence of 
perpetuum mobile). Controversies have also different sources, depending on 
the domain. In mathematics, disagreement comes mostly from different opinions 
on whether particular steps of a complex argument are correct, while in 
experimental sciences we should rather expect disagreements concerning such 
issues as the framework in which experiments were conducted,
statistical validity of experimental samples, or data interpretation.
Arguably, in the social sciences and humanities, the situation
may be very different, one of the reasons being that terms used in these 
domains are sometimes ambiguous 
or ill-defined, and as a result, even the crux of the controversy may be hard to specify. 
This issue, however, is perhaps better left for further discussion.

{\bf Acknowledgements.}
Thanks are due to Jody Azzouni for enlightening discussions concerning the subject of this paper and
to an anonymous referee whose important remarks permitted us to correct the arguments and improve presentation.
The remaining flaws are, of course, solely the responsibility of the author.

{\bf Funding.}
This research was partially supported by NSERC discovery grant 
and by the Research Chair in Distributed Computing at the
Universit\'e du Qu\'ebec en Outaouais.

\bigskip

{\bf \Large References}

Aschbacher, M. [2004]:
The status of the classification of the finite simple groups, 
Notices of the American Mathematical Society (7) 51, 736-740.

Azzouni, J. [2004]:  The derivation-indicator view of mathematical practice, 
Philosophia Mathematica (3) 12, 81-105.

Azzouni, J. [2006]: Tracking Reason. Proof, Consequence, and Truth, Oxford University Press.

De Bruijn, N.G. [1970]:  The mathematical language AUTOMATH, its usage and some of its extensions, 
Lecture Notes in Mathematics; 125. Berlin: Springer, 29-61.

De Millo, R.A., Lipton, R.J., and Perlis, A.J. [1979]: Social processes and proofs of theorems and programs, 
Communications of the ACM 22, 271-280.

Grabiner, J.V. [1974]: Is mathematical truth time-dependent?, American Mathematical Monthly 81, 354-365.

Halliday, D., Resnick, R., Walker, J. [1996]: Fundamentals of Physics. New York: Willey.

Hersh, R. [1979]:  Some proposals for reviving the philosophy of mathematics, Advances in Mathematics 31, 31-50.

Hilbert, D. [1900]: Mathematische Probleme, 
Nachrichten Konigl. Gesell. Wiss. Gottingen, Math.-Phys. Klasse, 253-297; 
English trans. (Mary W. Newson): Mathematical problems, Bull. Amer. Math. Soc. 8 (1902), 437-179.

MacKenzie, D. [2001]: Mechanizing proof: computing, risk, and trust. Cambridge, Massachusetts: MIT Press.

Rav, Y. [1999]: Why do we prove theorems?,  Philosophia Mathematica (3) 7, 5-41.

Rav, Y. [2007]: A critique of a formalist-mechanist version of the justification of arguments 
in mathematicians proof practices, Philosophia Mathematica (3) 15, 291-320.  

Thomas, R.S.D. [1990]: Inquiry into meaning and truth, Philosophia Mathematica (1) 5, 73-87.

Thurston, W.P. [1994]: Letter to the Editors, Scientific American 270, No. 1, 5.

Wiles, A. [1995]: Modular elliptic curves and Fermat's Last Theorem, Annals of Mathematics 141, 443-551.

\end{document}